\documentclass[12pt]{amsart}
\usepackage{graphicx}
\usepackage{amssymb}
\usepackage{epstopdf}
\usepackage{hyperref}

\usepackage{tikz}
\usetikzlibrary{arrows,shapes,positioning}
\usetikzlibrary{decorations.markings}
\tikzstyle arrowstyle=[scale=1]
\tikzstyle directed=[postaction={decorate,decoration={markings,
    mark=at position .65 with {\arrow[arrowstyle]{stealth}}}}]
\tikzstyle reverse directed=[postaction={decorate,decoration={markings,
    mark=at position .65 with {\arrowreversed[arrowstyle]{stealth};}}}]
\input diagxy

\newtheorem{innercustomgeneric}{\customgenericname}
\providecommand{\customgenericname}{}
\newcommand{\newcustomtheorem}[2]{%
  \newenvironment{#1}[1]
  {%
   \renewcommand\customgenericname{#2}%
   \renewcommand\theinnercustomgeneric{##1}%
   \innercustomgeneric
  }
  {\endinnercustomgeneric}
}

\newcustomtheorem{customthm}{Theorem}
\newcustomtheorem{customlemma}{Lemma}

\theoremstyle{definition}

\begin{document}

\title[(A)symmetric Coorbital Configurations]{Symmetry and Asymmetry in the 1+N Coorbital Problem}
\author {Yiyang Deng, Marshall Hampton, and Zhiqiang Wang}

\address{Marshall Hampton \\
Department of Mathematics and Statistics \\
University of Minnesota, Duluth\\
1117 University Drive\\
Duluth, MN, 55812, USA}

\address{Yiyang Deng \\College of Mathematics and Statistics \\ Chongqing Technology and Business University \\ Chongqing 400067 \\ China}

\address{Zhiqiang Wang \\ 
College of Mathematics and Statistics\\
Chongqing University \\
Chongqing 401331 \\
China}

\begin{abstract}
The relative equilibria of planar Newtonian $N$-body problem become coorbital around a central mass in the limit when all but one of the masses becomes zero.  We prove a variety of results about the coorbital relative equilibria, with an emphasis on the relation between symmetries of the configurations and symmetries in the masses, or lack thereof.   We prove that in the $N=4$, $N=6$, and $N=8$ Newtonian coorbital problems there exist symmetric relative equilibria with asymmetric positive masses.   This result can be generalized to other homogeneous potentials, and we conjecture similar results hold for larger even numbers of infinitesimal masses.  We prove that some equalities of the masses in the $1+4$ and $1+5$ coorbital problems imply symmetry of a class of convex relative equilibria.  We also prove there is at most one convex central configuration of the symmetric $1+5$ problem.
\end{abstract}

\maketitle

\section{Introduction}

The special case of one large mass and $N$ infinitesimal masses interacting through the Newtonian gravitational model (the `1+N'-body problem) was first formally investigated by Maxwell \cite{brush1983maxwell}, who was mainly interested in the case of very large $N$ to model planetary rings such as Saturn's.  More interest in the problem for small $N$ developed partially because of the discovery of the nearly coorbital satellites Janus and Epimetheus of Saturn \cite{yoder1983theory}, which travel on an interesting horseshoe orbit.  The preprint by Hall \cite{hall88} framed the problem more generally within the context of the study of central configurations, and around the same time Salo and Yoder \cite{salo1988dynamics} numerically extensively investigated the configurations and dynamics of the $1+N$ problem for $N<10$ for identical infinitesimal masses.    Verrier and McInnes have also analyzed tadpole and horsehoe orbits by numerical continuations from relative equilibria in the $1+2$, $1+3$, and $1+4$ coorbital problems \cite{verrier2014periodic}.  

Our observations in this paper extend the approach of Renner and Sicardy \cite{RS2004} in considering the antisymmetry of the mass coefficient matrix of the defining equations for these central configurations.  Similar considerations have been very useful in analyzing the collinear $N$-body central configurations as well \cite{albouy_inverse_2000}.

This manuscript focuses on the $1+N$ problem for $N \le 8$.   Numerically it appears that there is a single $1+N$ configuration when the $N$ masses are equal for $N \ge 9$ (the regular $N$-gon around the central mass), but lower values of $N$ exhibit more complexity.  Hall \cite{hall88} was able to prove the uniqueness of the $1+N$ equal mass central configuration for $N \ge e^{27000}$, which was improved to $N \ge e^{73}$ by Casasayas et al \cite{casasayas1994central}.  

% add in discussion of finite dominant mass

When the central mass is large but still finite relative to the coorbital small masses there are some additional results, all of which require equal small masses on a regular polygonal ring.  Maxwell showed that for a sufficiently large central mass this configuration would be stable\cite{brush1983maxwell}, although his analysis was incorrect for $N < 7$.  The stability for $N \ge 7$ was shown by Moeckel \cite{moeckel_dom94}.  Some related results on stability for regular polygonal $1+N$ configurations are proved in \cite{xu2013linear}.  

Very little is known about $1+N$ central configurations with unequal masses for $N \ge 4$.  

 The $1+N$ problem can be extended to other potentials, such as the point vortex model \cite{HelmholtzV, kirchhoff_vorlesungen_1883, barry10}.  There are some results on $1+3$  \cite{barry2016existence} and  $1+4$ \cite{hoyer2020symmetry,oliveira2020stability} planar vortex central configurations. The $1+3$ problem already possesses a remarkably rich structure, which has been well-studied in the Newtonian case \cite{corbera2011central, corbera2015bifurcation}.   The central configurations for the Newtonian $1+4$ problem with identical small masses were rigorously determined by Albouy and Fu \cite{albouy2009relative}.
 
\section{Equations for $1+N$ central configurations and their stability}

We consider it worthwhile to consider the $1+N$ central configuration problem in the more general case of a homogeneous potential.  The derivation of the equations in \cite{corbera2015bifurcation} extends easily to this setting.  Namely for a configuration of $N$ points to form a $1+N$ planar central configuration they must lie on a common circle, which can be scaled to be the unit circle.  Using polar coordinates with angles $\theta_i$, and using the shorthand $\theta_{i,j} = \theta_i - \theta_j$, we can write the equations in terms of

$$ f_{i,j} = \sin(\theta_{i,j})(\frac{1}{r_{i,j}^s} -1) $$ 
where $s$ is the potential exponent ($s=3$ is the Newtonian case, and $s=2$ the vortex case) and $r_{i,j}$ is the distance from point $i$ to point $j$.  We will abuse/overload our notation a little and also sometimes think of $f$ as a function of a single variable:

$$ f(\theta) = \sin(\theta)(2^{-s}|\sin{\theta/2}|^{-s} -1) $$

Note that for points on the unit circle there are a few equivalent ways to write the mutual distances:

$$ r_{i,j} = \sqrt{2 - 2 \cos(\theta_{i,j})} = 2 \left | \sin(\frac{\theta_{i,j}}{2}) \right | $$

The matrix $F$ with entries $f_{i,j}$ (with $f_{i,i} = 0$) is the mass coefficient matrix for the coorbital equations:

$$ F m = 0 $$

where $m = (m_1, m_2, \ldots, m_N)^T$. 

Since $F$ is real and antisymmetric, it has purely imaginary eigenvalues, and the dimension of the kernel of $F$ has the same parity as $N$.

We can view the coorbital relative equilibrium equations as conditions to have a critical point of the effective potential (sometimes referred to as Hall's potential): 

$$ V = \sum_{i<j} m_i m_j \left( \frac{1}{(s-2) r_{i,j}^{s-2}} + \frac{r_{i,j}^2}{2}\right) $$

restricted to points on the unit circle, since 

$$ \frac{\partial V}{\partial \theta_i} = -m_i \sum_{j \neq i} m_j f_{i,j} $$

This can be extended to the point vortex case \cite{barry10,barry2016existence} by using the potential

$$ V = \sum_{i < j} m_i m_j \left( \log(r_{i,j})+ \frac{r_{i,j}^2}{2}\right) $$

In a slight abuse of notation we will refer to this vortex potential as the $s=2$ case.

We define a diagonal mass matrix $M = diag(m_1,m_2, \ldots, m_N)$ so that the coorbital central configuration equations can be written

$$ F m = -M^{-1} \nabla V  = 0 $$

We will also denote the Hessian matrix of $V$ by $H$, i.e. with matrix entries

$$H_{i,j} = \frac{\partial^2 V}{\partial \theta_i \partial \theta_j} = m_i m_j \left ( -\cos(\theta_{i,j}) - \frac{s + (s-2)\cos(\theta_{i,j})}{2^{s+1}|\sin((\theta_{i,j})/2)|^{s}} \right ) $$

$$H_{i,i}=-\sum_{j\neq i}H_{i,j}$$

which can also be written in terms of the distances as

$$H_{i,j} = m_i m_j \left ( \frac{r_{i,j}^2}{2} - 1 + \frac{(s-2)}{4 r_{i,j}^{s-2}} - \frac{(s-1)}{r_{i,j}^s} \right ) = m_i m_j h_{i,j}$$

\subsection{Lemma on stability of $1+N$ central configurations}
 
Moeckel proved that for a sufficiently large central mass a relative equilibrium is linearly stable if and only if it is a local minimum of the potential $V$ \cite{moeckel_dom94}.  Renner and Sicardy \cite{RS2004} showed that a relative equilibrium of the coorbital problem is linearly stable if and only the matrix $A = M^{-1} H$ has nonpositive eigenvalues (with exactly one zero eigenvalue that arises from the rotational symmetry), and more generally each positive eigenvalue of $A$ implies a negative and a positive eigenvalue of the linearization around the relative equilibrium.  We extend their result with the following:

\begin{customlemma}{1}
If a relative equilibrium of the coorbital problem has Morse index $m$ (the largest dimension on which $H$ is negative definite), then the number of imaginary eigenvalues of the linearization of the equations of motion at the relative equilibrium is $2m$, and the number of positive and negative eigenvalues are both $N - m$.  
\end{customlemma}

We assume the masses are all positive, so that $M^{-1}$ has square root $$R = diag(1/\sqrt{m_1}, \ldots, 1/\sqrt{m_N})$$ and then

$$ A = M^{-1} H = R R H = R R H R R^{-1} = R (R^T H R) R^{-1} $$
and we see that $A$ is similar to $R^T H R$, so they have the same eigenvalues, and $R^T H R$ is congruent to $H$, so $A$ and $H$ have the same numbers of positive, zero, and negative eigenvalues.

\qed

\subsection{Some new equations for the 1+$N$ problem}

An alternative set of equations for the $1+N$ problem for $N>2$ can be obtained by using only mutual distances, and determining the constrained critical points of the potential $V$  through Lagrange multipliers.  If three points $a$, $b$, and $c$ are on the unit circle, then their mutual distances satisfy the equation

$$ E_{a,b,c} := (r_{a,b} r_{b,c} r_{a,c})^2 + r_{a,b}^{4} - 2 r_{a,b}^{2} r_{a,c}^{2} + r_{a,c}^{4} - 2 r_{a,b}^{2} r_{b,c}^{2} - 2 r_{a,c}^{2} r_{b,c}^{2} + r_{b,c}^{4} = 0 $$

For $N$ coorbital masses, there are $\binom{N}{2}$ mutual distances and $\binom{N}{3}$ 3-body constraints.  To be a critical point of the potential $V$ the configuration must also satisfy the $\binom{N}{2}$ critical point equations

$$ \frac{\partial V}{\partial r_{i,j}} + \sum_{k \neq i,j} \frac{\partial E_{i,j,k}}{\partial r_{i,j}}  \lambda_{i,j,k}  = 0 $$

with $\binom{N}{3}$ Lagrange multipliers $\lambda_{i,j,k}$.  

For $N=3$ these equations are fairly simple, since there is only one constraint $E_{1,2,3} = 0$.  The three critical point equations (with $\lambda = \lambda_{1,2,3}$), in terms of $W_{i,j} = m_i m_j (1 - \frac{1}{r_{i,j}^{s}}) \}$, are

$$ r_{1,2} \left ( W_{1,2} + 2 \lambda   (r_{1,3}^{2} r_{2,3}^{2} + 2 r_{1,2}^{2} - 2 r_{1,3}^{2} - 2 r_{2,3}^{2})  \right ) = 0 $$

$$ r_{1,3} \left ( W_{1,3} + 2 \lambda (r_{1,2}^{2} r_{2,3}^{2} + 2 r_{1,3}^{2} - 2 r_{1,2}^{2} - 2 r_{2,3}^{2})    \right ) = 0 $$

$$ r_{2,3} \left ( W_{2,3} + 2 \lambda  (r_{1,2}^{2} r_{1,3}^{2} + 2 r_{2,3}^{2} - 2 r_{1,2}^{2} - 2 r_{1,3}^{2})   \right )  = 0 $$

and we could in addition divide each equation by the $r_{i,j}$ factor.

A second approach to obtain rational or polynomial equations in the mutual distances is to choose a set of independent distances and use geometric relations to compute the partial derivatives of the remaining distances with respect to the basis set.  Let $W_{i,j} = m_i m_j (1 - \frac{1}{r_{i,j}^{s}}) \}$, and then

$$ \frac{\partial V}{\partial r_{1,i}} = r_{1,i} W_{1,i}  + \sum_{j}  r_{i,j} W_{i,j} \frac{\partial r_{i,j}}{\partial r_{i,i+1}} =  0$$
and then we can divide out the $r_{1,i}$.

This approach becomes unwieldy for larger $N$, but can be useful for particular cases in which there are additional restrictions on the configuration.  As an example, for $N=3$, if our set of independent distances is $r_{1,2}$ and $r_{2,3}$, we can compute $\frac{\partial r_{1,3}}{\partial r_{1,2}}$ and $\frac{\partial r_{1,3}}{\partial r_{2,3}}$  from $E_{1,2,3}=0$. The equations then become

$$ r_{1,2} W_{1,2} - r_{1,3} W_{1,3} \frac{\frac{\partial E_{1,2,3}}{\partial r_{1,2}}}{\frac{\partial E_{1,2,3}}{\partial r_{1,3}}} = r_{1,2} \left ( W_{1,2} - W_{1,3} \frac{{\left(r_{1,3}^{2} r_{2,3}^{2} - 2 \, r_{1,2}^{2} + 2 \, r_{1,3}^{2} + 2 \, r_{2,3}^{2}\right)} }{{\left(r_{1,2}^{2} r_{2,3}^{2} + 2 \, r_{1,2}^{2} - 2 \, r_{1,3}^{2} + 2 \, r_{2,3}^{2}\right)} }  \right ) = 0$$

$$ r_{2,3} W_{2,3} - r_{1,3} W_{1,3} \frac{\frac{\partial E_{1,2,3}}{\partial r_{1,2}}}{\frac{\partial E_{1,2,3}}{\partial r_{2,3}}} = r_{2,3} \left ( W_{2,3} - W_{1,3} \frac{{\left( r_{1,3}^{2} r_{2,3}^{2} - 2 \, r_{1,2}^{2} + 2 \, r_{1,3}^{2} + 2 \, r_{2,3}^{2} \right)} }{{\left(r_{1,2}^{2} r_{1,3}^{2} + 2 \, r_{1,2}^{2} + 2 \, r_{1,3}^{2} - 2 \, r_{2,3}^{2}\right)} } \right )
 = 0$$
 
For the $N=3$ case these two approaches are essentially the same, as it is simple to eliminate $\lambda$ from the first set of equations, but they diverge more significantly for larger $N$.

\section{The 1+2(+1) Problem}

In the case of two infinitesimal masses it is elementary to determine that the relative equilibria are determined by the zeros of $f_{1,2}$, with $\theta_{1,2} \in \{ \pi/3, \pi, 5\pi/3\}$, independent of the masses.  These correspond to the non-coalescent zero-mass limits of the well-known Euler and Lagrange central configurations.  

%add Euler and Lagrange references

% Stability? - as special cases of Euler and Lagrange

Now we consider adding an additional infinitesimal particle to the 1+2 setting so that the resulting `1+2+1' configuration is still a relative equilibrium.  Let the angle of this particle be $\theta_3$, and normalize the masses of the initial two masses by $m_1 + m_2 = 1$.  Then $\theta_3$ must satisfy

$$  m_1 f_{1,3} + (1-m_1) f_{2,3} = 0 $$

\begin{customlemma}{2}
The 1+2+1 problem has ten solutions for any Hall potential with $s>2$ and positive masses $m_1$ and $m_2$.
\end{customlemma}
\begin{proof}
Without loss of generality we can set $\theta_1 = 0$ and then either $\theta_2 = \pi/3$ ($\theta_2 = 5\pi/3$ can be viewed as a relabeling of the $\theta_2 = \pi/3$ case) or $\theta_2 = \pi$.

In the first case, $\theta_2 = \pi/3$, the equation for $\theta_3$ becomes

$$ g_1 = m_1 \sin(\theta_3) \left (\frac{1}{2^s | \sin(\theta_3/2)|^s} - 1 \right ) + (1-m_1) \sin(\theta_3-\pi/3) \left (\frac{1}{2^s | \sin(\theta_3/2 - \pi/6)|^s} - 1 \right )  = 0 $$

An elementary calculation shows the derivative of $g_1$ is negative on the interval $(0, \pi/3)$ for $m_1 \in (0,1)$. At the endpoints,  $g_1$ is positive for sufficiently small $\theta_3 > 0$, and is negative for $\theta_3  = \pi/3 - \epsilon$ for sufficiently small $\epsilon > 0$.  So for each $s>2$ and $m_1 \in (0,1)$ there is a unique solution $\theta_3 \in (0, \pi/3)$.

In the interval $(\pi/3, \pi)$ the second derivative of $g_1$ is positive for $s>2$ and $m_1 \in (0,1)$.  The derivative of $g_1$ changes sign on this interval, so there is a unique zero for $d g_1/ d \theta_3$ in $(\pi/3, \pi)$.  It is also easy to check that for sufficiently small $\epsilon > 0$ the sign of $g_1$ is positive for $\pi/3 + \epsilon$ and negative for $\pi - \epsilon$, so there is also a unique zero of $g_1$ for $\theta_3 \in (\pi/3, \pi)$.

For $\theta_3$ in the interval $(\pi, 4 \pi/3)$, the function $g_1$ has a change of sign, and its derivative is positive, so the solution in this interval is unique.  

The final interval, $(4 \pi/3, 2 \pi)$, has the same properties as the interval $(\pi/3, \pi)$ (after interchanging the masses, and an overall sign change from reversing the angular differences) and it contains another unique solution.

In the second case, with $\theta_2 = \pi$, the angle $\theta_3$ must satisfy

$$ g_1 = \sin(\theta_3) \left [ m_1 \left (\frac{1}{2^s | \sin(\theta_3/2)|^s} - 1 \right ) - (1-m_1)  \left (\frac{1}{2^s | \sin(\theta_3/2 - \pi/2)|^s} - 1 \right ) \right ] = 0 $$

The intervals $(0,\pi)$ and $(\pi,2 \pi)$ have symmetric behavior in this problem, so we need only consider $\theta_3 \in (0,\pi) $.  Since $\sin(\theta_3)$ is not zero on $(0,\pi)$, we can simplify the condition to

$$ \tilde{g}_1 =  m_1 \left (\frac{1}{2^s | \sin(\theta_3/2)|^s} - 1 \right ) - (1-m_1)  \left (\frac{1}{2^s | \sin(\theta_3/2 - \pi/2)|^s} - 1 \right ) = 0 $$

The function $\tilde{g}_1$ is monotonically decreasing on $(0,\pi)$ and has a change of sign on that interval, so there is a unique solution in $(0,\pi)$.  By the symmetry of the problem, there is also a unique solution in $(\pi,2 \pi)$.  

Altogether there are 8 solutions, 3 with $\theta_2 = \pi/3$, three with $\theta_2 = 5 \pi/3$, and two with $\theta_2 = \pi$.

\end{proof}

\section{The 1+2N Problem}
\subsection{The 1+4 Problem}

In the $1+4$ problem complementary observations to ours on symmetrical configurations have been shown by Oliveira \cite{oliveira2013symmetry}, who analyzed the case in which two masses are at opposite points of the circle.  In this case, to be a relative equilibrium the configuration must be symmetric, and have symmetric masses (i.e. if masses 1 and 3 are at opposite points of the circle, then their axis is an axis of symmetry  and $m_2 = m_4$).    We will call these type-1 symmetric configurations of the 1+4 problem.  Without loss of generality we can assume that $\theta_1 = 0$, $\theta_3 = \pi$, and $\theta_2 = -\theta_4 \in (0,\pi)$.

\begin{customlemma}{3}
Every type-1 symmetric configuration of the 1+4 problem has a two dimensional space of real mass solutions.
\end{customlemma}
This follows immediately from the structure of the mass coefficient matrix, which has the form

$$ \left(\begin{array}{rrrr}
0 & {f_{1,2}} & 0 & -{f_{1,2}} \\
-{f_{1,2}} & 0 & {f_{2,3}} & {f_{2,4}} \\
0 & -{f_{2,3}} & 0 & {f_{2,3}} \\
{f_{1,2}} & -{f_{2,4}} & -{f_{2,3}} & 0
\end{array}\right) $$

Any matrix of this form has a vanishing determinant.  Since $f_{1,2} \neq 0$ and $f_{2,3} \neq 0$, it has at least rank 1.  As an antisymmetric real matrix the kernel must be even dimensional.  Together these imply that the kernel must be two-dimensional.  It is also not difficult to compute an explicit basis for this kernel, for example the mass vectors $(f_{2,3},0,f_{1,2},0)$ and $(0,f_{2,3},-f_{2,4},f_{2,3})$.  

Quite surprisingly we obtain a very different result for the other type of symmetric 1+4 configurations.  We will define type-2 symmetric configurations of the 1+4 problem as having a symmetry axis through the center of the circle with two pairs of points off of the symmetry axis.  For representatives of these configurations we will assume $\theta_4 = -\theta_1$ and $\theta_3 = -\theta_2$, so $f_{3,4} = f_{1,2}$ and $f_{2,4} = f_{1,3}$.

Some closely related results were obtained by Deng, Li, and Zhang \cite{deng2019symmetric}, who focused on the forward problem of determining the number of 1+4 relative equilibria given the masses, rather than the inverse problem of determining the masses from the configuration.  

% Put in conjecture of Deng,Li, Zhang somewhere

\begin{customthm}{1} \label{n4}
For the $N=4$ Newtonian coorbital problem there exist type-2 symmetric central configurations with asymmetric positive masses. 
\end{customthm}

\begin{figure}[h!t]
\includegraphics[width=3in]{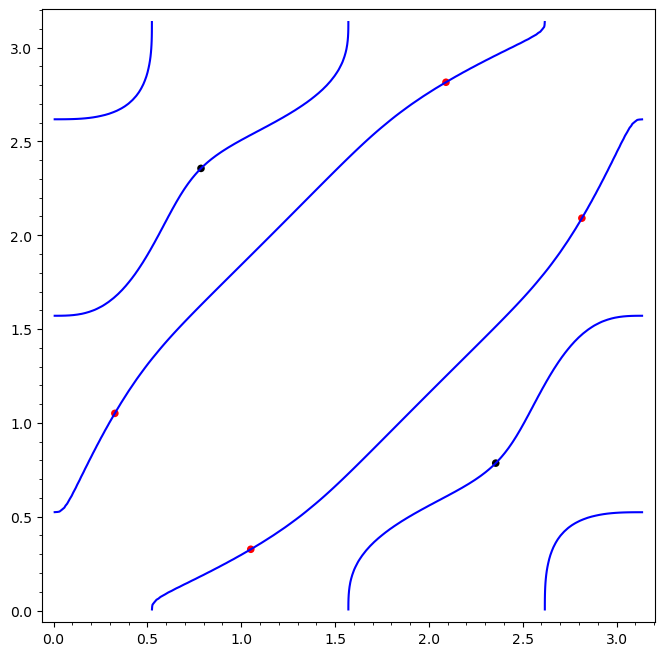}
\caption{The zero set of the Pfaffian for our symmetric 1+4 configurations.  The locations of the square configuration are marked by black dots and the convex equal mass configurations by red dots.}
\label{pfaff4}
\end{figure}

The matrix $F$ for $N=4$ has the Pfaffian $ f_{1,2} f_{3,4} + f_{1,4} f_{2,3} - f_{1,3} f_{2,4}$.    For our type-2 symmetric configurations the Pfaffian simplifies to $f_{1,2}^2 - f_{1,3}^2 + f_{1,4}f_{2,3}$.  In Figure \ref{pfaff4} we show a numerical computation of the zero locus of this Pfaffian for the Newtonian potential.

Numerically we find that there is a one-dimensional subset of angles $(\theta_1, \theta_2)$ with asymmetric masses but for this proof we will simply show one example.

If $\theta_1 = \frac{\pi}{6}$, then $r_{1,4} = 1$ and $f_{1,4} = 0$.  For the Pfaffian to vanish, we need $f_{1,2}^2 - f_{1,3}^2 = 0$.  It is easy to verify that there is an angle $\theta_2 \approx 1.936$ such that $f_{1,2} = - f_{1,3} \approx -0.536$ (with interval arithmetic for example, or by bracketing this value with exactly computable quantities).  

For these angles the matrix $F$ then has the form

$$ \left(\begin{array}{rrrr}
0 & f_{1,2} & -f_{1,2} & 0 \\
-f_{1,2} & 0 & f_{2,3} & -f_{1,2} \\
f_{1,2} & - f_{2,3} & 0 & f_{1,2} \\
0 & f_{1,2} & -f_{1,2} & 0
\end{array}\right) $$

A basis for the two-dimensional kernel is $v_1 = (1, 0, 0, -1)$, $v_2 =(0, f_{1,2}, f_{1,2}, f_{2,3})$.  It is again easy to check that $f_{2,3}$ is negative for these angles ($f_{2,3} \approx -0.565$), so for $\alpha \in (-\infty, 1/f_{2,3})$ we have positive mass solution vectors $v_1 + \alpha v_2$.  These are asymmetric ($m_1 \neq m_4$) for all $\alpha \neq 2/f_{2,3}$.

\begin{center}
\begin{figure}[h!t]
\includegraphics[width=2.0in]{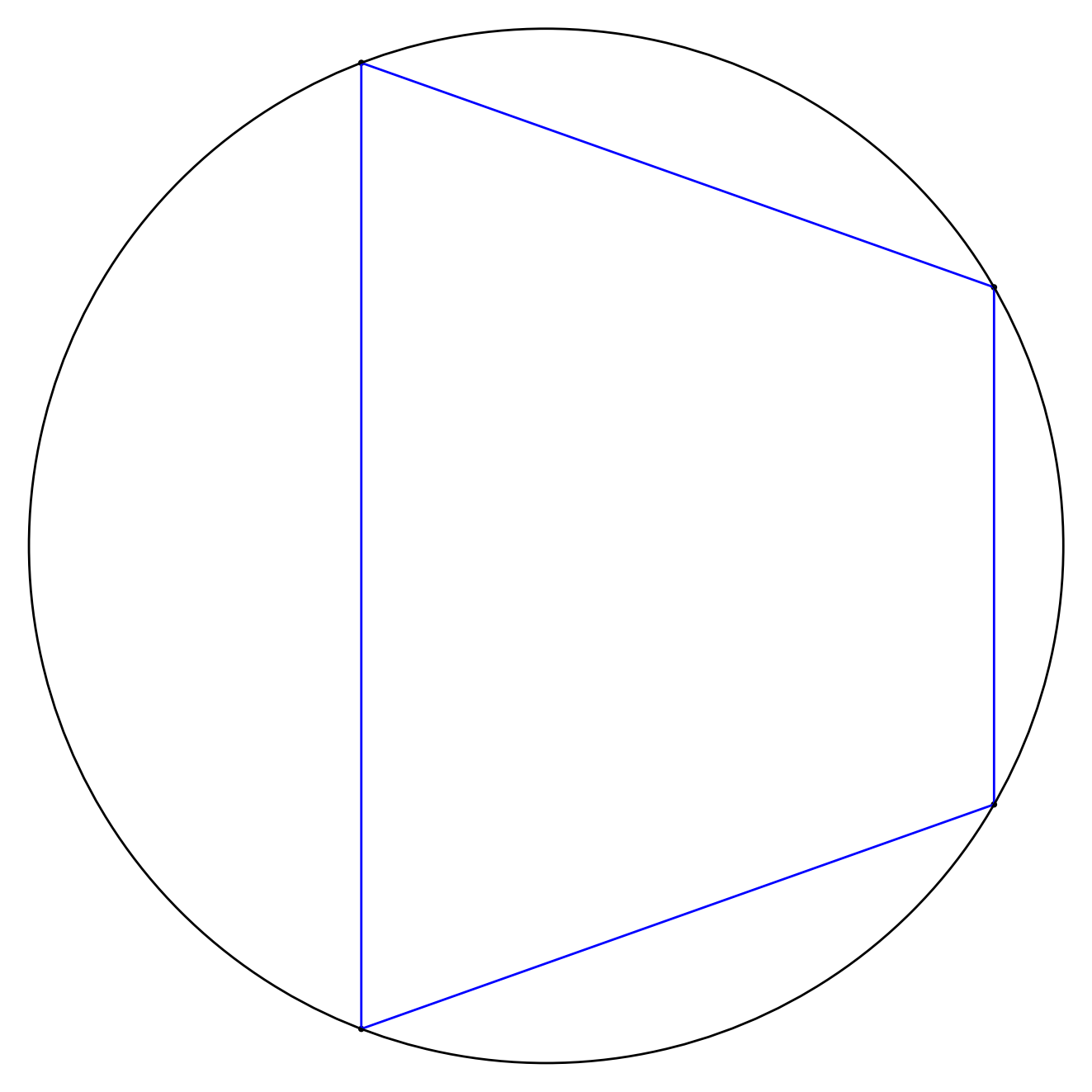}
\caption{An example of a symmetric 1+4 central configuration with asymmetric masses}
\end{figure}
\end{center}

%%%% Begin Yiyang theorems for 1+ 4

In contrast to the above result, in some cases if there is a symmetry in the masses there must a corresponding symmetry in the configuration.  The following lemma is a key tool in our proofs of these results, providing a factorization of a commonly occurring term from the central configuration equations.

\begin{customlemma}{4} \label{c4lem}
For a convex four-point coorbital configuration with $\theta_{i,j}$, $\theta_{k,l}$ and $\theta_{i,l}$ in $[0,\pi/2]$, in the Newtonian case ($s=3$) the quantity 

$$ r_{i,j}^3 r_{k,l}^3 ( f_{i,j} -  f_{k,l} ) = \sin(\theta_{i,j}) r_{k,l}^3 (1 - r_{i,j}^3) - \sin(\theta_{k,l}) r_{i,j}^3 (1 - r_{k,l}^3)$$

 can be written as 

\begin{align*}
& 32 \sin (\frac{\theta_{i,j}}{2}) \sin (\frac{\theta_{k,l}}{2}) \sin \frac{(\theta_{k,l}  - \theta_{i,j})}{4}  \\
 & \left[  8 \sin^2 (\frac{\theta_{i,j}}{2}) \sin^2 (\frac{\theta_{k,l}}{2})  \cos\frac{(\theta_{i,j} + \theta_{k,l})}{2} \cos  \frac{(\theta_{k,l}  - \theta_{i,j})}{4}  + \sin  \frac{(\theta_{k,l}  + \theta_{i,j})}{4} (1 + \cos (\frac{\theta_{i,j}}{2}) \cos (\frac{\theta_{k,l}}{2}) ) \right] \\
 % & =: \sin \frac{(\theta_{k,l}  - \theta_{i,j})}{4}  G_{i,j;k,l}
\end{align*}

% and the quantity $G_{i,j;k,l}$ is positive.

\end{customlemma}

\begin{proof}

With the assumptions on the angles we can write the distances in terms of trigonometric functions and factor out some common terms:

\begin{align*}
       &\sin(\theta_{i,j}) r_{k,l}^3 (1 - r_{i,j}^3) - \sin(\theta_{k,l}) r_{i,j}^3 (1 - r_{k,l}^3) \\
     = & 16 (1 - 8\sin^3 \frac{\theta_{i,j}}{2}) \sin^3 \frac{\theta_{k,l}}{2} \sin \frac{\theta_{i,j}}{2} \cos \frac{\theta_{i,j}}{2}  - 16 (1 - 8\sin^3 \frac{\theta_{k,l}}{2}) \sin^3 \frac{\theta_{i,j}}{2} \sin \frac{\theta_{k,l}}{2} \cos \frac{\theta_{k,l}}{2} \\
     = & 16 \sin \frac{\theta_{i,j}}{2} \sin \frac{\theta_{k,l}}{2} \bigg[(1 - 8 \sin^3 \frac{\theta_{i,j}}{2}) \sin^2 \frac{\theta_{k,l}}{2} \cos
         \frac{\theta_{i,j}}{2}  -(1 - 8 \sin^3 \frac{\theta_{k,l}}{2}) \sin^2 \frac{\theta_{i,j}}{2} \cos \frac{\theta_{k,l}}{2} \bigg].
\end{align*}

and now we continue with some elementary trigonometric identities (sum-to-product) and refactoring:

\begin{align*}
       =  & 16 \sin (\frac{\theta_{i,j}}{2}) \sin (\frac{\theta_{k,l}}{2}) \left[ 8 \sin^2 (\frac{\theta_{i,j}}{2}) \sin^2 (\frac{\theta_{k,l}}{2}) \cos \frac{( \theta_{i,j} + \theta_{k,l})}{2} \sin \frac{( \theta_{k,l} - \theta_{i,j})}{2} \right . \\
       & \left . - (\cos (\frac{\theta_{k,l}}{2}) - \cos (\frac{\theta_{i,j}}{2})) (1 + \cos (\frac{\theta_{i,j}}{2}) \cos (\frac{\theta_{k,l}}{2})) \right] \\
      = & 32 \sin (\frac{\theta_{i,j}}{2}) \sin (\frac{\theta_{k,l}}{2}) \sin \frac{( \theta_{k,l} - \theta_{i,j})}{4} \left[8 \sin^2 (\frac{\theta_{i,j}}{2}) \sin^2 (\frac{\theta_{k,l}}{2}) \cos \frac{( \theta_{i,j} + \theta_{k,l})}{2} \cos\frac{( \theta_{k,l} - \theta_{i,j})}{4} \right . \\
      & \left . + \sin \frac{( \theta_{i,j} + \theta_{k,l})}{4} (1 + \cos (\frac{\theta_{i,j}}{2}) \cos (\frac{\theta_{k,l}}{2})) \right]
\end{align*}

% The positivity of the quantity $G_{i,j;k,l}$ directly follows from the assumptions on the angles.
\end{proof}

We will use this lemma in the following theorem (and later for a similar result in the 1+5 symmetric case):

\begin{customthm}{2} \label{sym4}
A 1+4 convex coorbital central configuration with $m_1 = m_4>0$ and $m_2 = m_3>0$ must be symmetric.
\end{customthm}
\begin{proof} 

Without loss of generality we can assume that the angles of the convex configuration satisfy $-\frac{\pi}{2} < \theta_4 < \theta_3 < 0 < \theta_2 < \theta_1,$ and $\theta_{1,4} \leq \pi$.

We can sum the second and third rows of the central configuration equations ($Fm = 0$), assuming $m_1 = m_4$ and $m_2 = m_3$, to cancel out the $m_2$ terms and obtain

\begin{equation}\label{eq21} \sin(\theta_{1,2}) (\frac{1}{r_{1,2}^3} - 1) + \sin(\theta_{4,2}) (\frac{1}{r_{2,4}^3} - 1) + \sin(\theta_{1,3}) (\frac{1}{r_{1,3}^3} - 1) + \sin(\theta_{4,3}) (\frac{1}{r_{3,4}^3} - 1) = 0. 
\end{equation}

% A1 + B1 = \frac{(\theta_{3,4} + \theta_{1,2})}{2}
% A2 + B2 = \frac{(\theta_{2,4} + \theta_{1,3})}{2}

Using Lemma 4 we can rewrite the inner and outer pairs of terms in the above relation, and then use the fact that $\frac{(\theta_{2,4} - \theta_{1,3})}{2} = \frac{(\theta_{3,4} - \theta_{1,2})}{2} = \frac{(\theta_2 + \theta_3) - (\theta_1 + \theta_4)}{2}$ to rewrite the equation (\ref{eq21}) as:
\begin{equation*}
  \begin{split}
   & \sin \frac{(\theta_{3,4} - \theta_{1,2})}{4} \Big[ 4\cos \frac{(\theta_{3,4} - \theta_{1,2})}{4} \big(\cos \frac{(\theta_{3,4} + \theta_{1,2})}{2} + \cos \frac{(\theta_{2,4} + \theta_{1,3})}{2} \big) + \\
   & \frac{1}{2 \sin^2 \frac{\theta_{1,2}}{2} \sin^2 \frac{\theta_{3,4}}{2}} \sin \frac{(\theta_{3,4} + \theta_{1,2})}{4}(1 + \cos \frac{\theta_{1,2}}{2} \cos \frac{\theta_{3,4}}{2}) + \\
   & \frac{1}{2 \sin^2 \frac{\theta_{1,3}}{2} \sin^2 \frac{\theta_{2,4}}{2}} \sin \frac{(\theta_{2,4} + \theta_{1,3})}{4}(1 + \cos \frac{\theta_{1,3}}{2} \cos \frac{\theta_{2,4}}{2})\Big] =0.
  \end{split}
\end{equation*}

It is elementary to check that the quantities $\sin \frac{\theta_{1,2}}{2}$, $\sin \frac{\theta_{3,4}}{2}$,  $\cos \frac{(\theta_{3,4} - \theta_{1,2})}{4}$, $\sin \frac{(\theta_{1,2} + \theta_{3,4})}{4}$, $\sin \frac{\theta_{1,3}}{2}$, $\sin \frac{\theta_{2,4}}{2}$, $\cos \frac{(\theta_{2,4} - \theta_{1,3})}{4}$, and $\sin \frac{(\theta_{1,3} + \theta_{2,4})}{4}$ are all positive for our assumptions on the angles.
 Since
\begin{align*}
   & \cos(\frac{\theta_{1,2}}{2} + \frac{\theta_{3,4}}{2}) + \cos(\frac{\theta_{1,3}}{2} + \frac{\theta_{2,4}}{2}) \\
  = & \cos(\frac{\theta_{1,4}}{2} - \frac{\theta_{2,3}}{2}) + \cos(\frac{\theta_{1,4}}{2} + \frac{\theta_{2,3}}{2}) \\
  = & \ 2 \cos\frac{\theta_{1,4}}{2} \cos\frac{\theta_{2,3}}{2} > 0,
\end{align*}
the terms included in the bracket are all positive. This implies that 
$$\sin \frac{\theta_{3,4} - \theta_{1,2}}{4}= 0,$$
which in turn implies that $\frac{\theta_{3,4}}{2} = \frac{\theta_{1,2}}{2}$. We conclude that the configuration must be symmetric.
\end{proof}

For positive mass solutions we must have some differences in sign in the mass coefficient in the coorbital central configuration equations. For analyzing the convex configurations (with $-\frac{\pi}{2} < \theta_4 < \theta_3 < 0 < \theta_2 < \theta_1,$ and $\theta_{1,4} \leq \pi$), we get the following simple bound:
\begin{customlemma}{5} The convex symmetric 1+4 central configurations with $-\frac{\pi}{2} < \theta_4 < \theta_3 < 0 < \theta_2 < \theta_1,$ and $\theta_{1,4} \leq \pi$ are contained in the set $\mathcal{B}$ defined by $\theta_{1,2} < \pi/3$, $\theta_{2,3} < \pi/3$, $\theta_{3,4} < \pi/3$ and $\theta_{1,4} > \pi/3$.

\end{customlemma}
\begin{proof} From the 1+4 coorbital central configuration equations $Fm = 0$, under our ordering of the masses we have $r_{1,2} < r_{1,3} < r_{1,4} \leq 2$, $r_{2,3} < r_{1,3}, r_{2,4}$ and $r_{43} < r_{42} < r_{41} \leq 2$. For the convex case, we can assume as above that $-\frac{\pi}{2} < \theta_4 < \theta_3 < 0 < \theta_2 < \theta_1,$ and $\theta_{1,4} \leq \pi$.

First we claim that $r_{1,2} < 1$. Indeed, if $r_{1,2}\geq 1$, then $1 \leq r_{1,2} < r_{1,3} < r_{1,4}$, and the left-hand side of the first row of  $Fm = 0$ would be negative. So $r_{1,2} < 1$, i.e $\sin (\frac{\theta_{1,2}}{2}) < \frac{1}{2}$ which implies $\theta_{1,2} < \pi/3$. Similarly, we also have $r_{3,4} < 1$ and $\theta_{3,4} < \pi/3$, and $r_{1,4} > 1$ and $\theta_{1,4} > \pi/3$.

Next we can show that $r_{2,3} < 1$. If we assume $r_{2,3} \geq 1$, then $1 \leq r_{2,3} < r_{2,4}$ and then the second row of  $Fm = 0$ would be positive.  So $r_{2,3} < 1$, which is equivalent to $\sin (\frac{\theta_{2,3}}{2}) < \frac{1}{2}$, which implies $\theta_{2,3} < \pi/3$.
\end{proof}

Furthermore, we get the same result as Theorem 1 in \cite{oliveira2013symmetry}:
\begin{customthm}{3} Under our assumption, if $m_1$ and $m_4$ are collinear with the central mass, then there is no coorbital central configuration for any given positive masses $m_1, m_2, m_3, m_4$.
\end{customthm}
\begin{proof} If $m_1$ and $m_4$ are collinear, then $\theta_{1,4} = \pi$. Since $\theta_{1,2} < \pi/3$ and $\theta_{3,4} < \pi/3$, we obtain that $\theta_{2,3} > \pi/3$, contradicting the previous lemma.
\end{proof}

%%%% END Yiyang for 1+4

\subsection{The 1+6 and 1+8 Problems}

For the symmetric configurations of the 1+6 problem we can again consider type-1 cases in which $\theta_1 = 0$ and $\theta_4 = \pi$, and $\theta_2 = -\theta_6$, $\theta_3 = -\theta_5$, and type-2 cases with 3 pairs of reflected points $\theta_1 = -\theta_6$, $\theta_2 = -\theta_5$, and $\theta_3 = -\theta_4$.

Similar to the 1+4 problem, the type-1 configurations have a mass coefficient matrix of the form

$$ \left(\begin{array}{rrrrrr}
0 & f_{1,2} & f_{1,3} & 0 & -f_{1,3} & -f_{1,2} \\
-f_{1,2} & 0 & f_{2,3} & f_{2,4} & f_{2,5} & f_{2,6} \\
-f_{1,3} & -f_{2,3} & 0 & f_{3,4} & f_{3,5} & f_{2,5} \\
0 & -f_{2,4} & -f_{3,4} & 0 & f_{3,4} & f_{2,4} \\
f_{1,3} & -f_{2,5} & -f_{3,5} & -f_{3,4} & 0 & f_{2,3} \\
f_{1,2} & -f_{2,6} & -f_{2,5} & -f_{2,4} & -f_{2,3} & 0
\end{array}\right) $$

which always has a zero determinant, so every such configuration has at least a 2-dimensional space of mass vector solutions (not necessarily positive).

\begin{customthm}{4}
For the $N=6$ and $N=8$ Newtonian coorbital problem there exist symmetric central configurations with asymmetric positive masses. 
\end{customthm}
The proof follows in a very similar way to Theorem \ref{n4}.  For $N=6$ we can use the particular case of $\theta_1 = \frac{\pi}{8}$, $\theta_2	=\frac{3 \pi}{7}$.   The Pfaffian changes sign for a $\theta_3 \approx 2.5349$ in $(\frac{4 \pi}{5}, \frac{13 \pi}{16})$, and one of the mass vector solutions can be chosen to be asymmetric. 

For $N=8$, we found a particular case contained in the intervals
$$ \frac{61 \pi }{702} < \theta_1 < \frac{63 \pi }{725} $$
$$ \frac{95 \pi }{289} < \theta_2 < \frac{24 \pi }{73} $$
$$ \frac{32 \pi }{55} < \theta_3 < \frac{71 \pi }{122} $$
$$ \frac{91 \pi }{110} < \theta_4 < \frac{24 \pi }{29} $$
for which the Pfaffian vanishes and the mass vector solutions can be chosen to be asymmetric. 

We used an ad-hoc approach to find these cases; a more systematic method might generalize to any even $N$.

\section{The 1+N Problem for odd N}

When $N$ is odd, the matrix $F$ must have a kernel of odd dimension, and so every configuration has a real nonzero (but not necessarily positive) vector of masses satisfying the relative equilibria conditions.  

For symmetric 1+5 configurations, for which we will choose representatives with $\theta_1 = - \theta_5$, $\theta_2 = -\theta_4$, and $\theta_3 = 0$, the matrix $F$ becomes

$$F = \left(\begin{array}{rrrrr}
0 & {f_{1,2}} & {f_{1,3}} & {f_{1,4}} & {f_{1,5}} \\
-{f_{1,2}} & 0 & {f_{2,3}} & {f_{2,4}} & {f_{1,4}} \\
-{f_{1,3}} & -{f_{2,3}} & 0 & {f_{2,3}} & -{f_{1,3}} \\
-{f_{1,4}} & -{f_{2,4}} & -{f_{2,3}} & 0 & {f_{1,2}} \\
-{f_{1,5}} & -{f_{1,4}} & {f_{1,3}} & -{f_{1,2}} & 0
\end{array}\right)$$

If we assume that the masses are also symmetric, with $m_1 = m_5$ and $m_2 = m_4$, then there are only two independent equations:

\begin{equation}\label{zmat}
 \left(\begin{array}{rrr}
{f_{1,5}} & {f_{1,2}} + {f_{1,4}} & {f_{1,3}} \\
-{f_{1,2}} + {f_{1,4}} & {f_{2,4}} & {f_{2,3}}
\end{array}\right)  \left(\begin{array}{r} m_1 \\ m_2 \\ m_3 \end{array}\right)  = 0 
\end{equation}

%To determine where we can have positive mass solutions for this symmetric case, we study the boundary where some masses are zero.  

\begin{center}
\begin{figure}[h!t]
\includegraphics[width=3in]{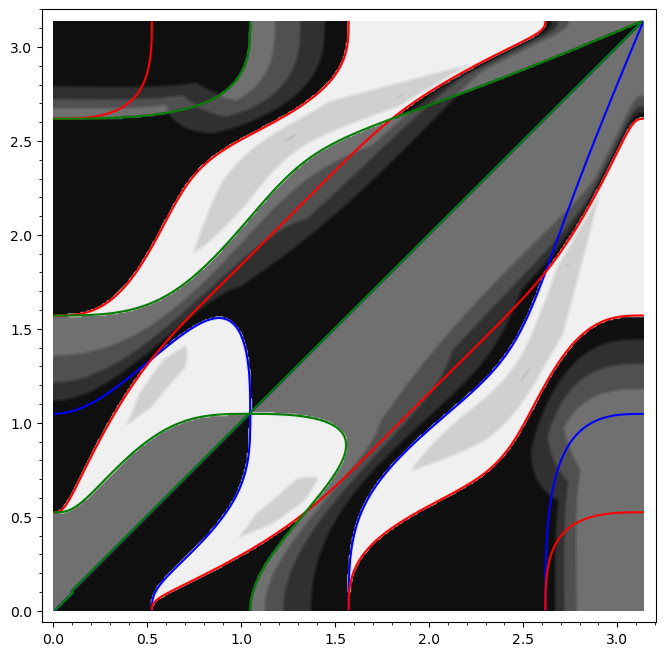}
\caption{Zero mass curves in the ($\theta_1$, $\theta_2$) plane for symmetric 1+5 configurations.  The red curve is where the central mass $m_3 = 0$, green and blue curves are $m_1=0$ and $m_2 = 0$.  Lighter colors correspond to large positive mass ratios, dark colors to large negative mass ratios.}
\label{zmfig}
\end{figure}
\end{center}

For a zero-mass configuration there must be a vanishing minor of the above coefficient matrix, i.e.

$$ Z_1 = {f_{2,3}}( {f_{1,2}} + {f_{1,4}}) - {f_{1,3}} {f_{2,4}}  = 0 $$
$$ Z_2 = {f_{1,3}}( {f_{1,2}}  - {f_{1,4}}) + {f_{1,5}} {f_{2,3}} = 0 $$
$$ Z_3 = {f_{1,2}}^{2} - {f_{1,4}}^{2} + {f_{1,5}} {f_{2,4}} = 0 $$

where $Z_i$ is must be zero if $m_i$ is zero.  The geometry of these curves is shown in Figure \ref{zmfig}.

Without loss of generality we can assume that $\theta_1 > \theta_2 > 0$.  For positive mass solutions we must have some differences in sign in the mass coefficient matrix above.  For analyzing the convex configurations (with $\theta_1 \le \pi/2$) we use the following simple bound:

\begin{customlemma}{6} 
The convex symmetric $1+5$ central configurations with $0 < \theta_2 < \theta_1 <  \pi/2$ are contained in the set $\mathcal{C}$ defined by $\pi/6  < \theta_2 < \pi/3$ and $\theta_{1,2} < \pi/3 $.
\end{customlemma}
\begin{proof}
With the assumption that $0 < \theta_2 < \theta_1 <  \pi/2$ the functions $f_{i,j}$ in the matrix \ref{zmat} are negative for $\theta_{i,j} < \pi/3$.  If $\theta_1 \le \pi/6$, then $\theta_2 < \pi/6$ as well, and then all of the entries of the first row of the mass coefficient matrix would be non-positive with at least one negative entry, so there could not be a mass vector in the kernel with all positive entries.  Similarly, if $\theta_2 \ge \pi/3$ then $\theta_1 > \pi/3$ and all of the entries of the first row would be non-negative with at least one positive entry.

\end{proof}

With this lemma we can prove the following

\begin{customthm}{5} 
There is at most one convex central configuration of the symmetric $1+5$ problem.
\end{customthm}

\begin{proof}

In the symmetric case, the Hessian $H$ is congruent to a block diagonal matrix using the matrix

$$P = \left(\begin{array}{rrrrr}
1 & 1 & 0 & 1 & 0 \\
1 & 0 & 1 & 0 & 1 \\
1 & 0 & 0 & 0 & 0 \\
1 & 0 & 1 & 0 & -1 \\
1 & 1 & 0 & -1 & 0
\end{array}\right) $$

so that

$$ P^T H P = \left(\begin{array}{rrrrr}
0 & 0 & 0 \\
0 & \mathcal{H}_1 & 0 \\
0 & 0 & \mathcal{H}_2 \\
\end{array}\right) $$

in which 

$$ \mathcal{H}_1 =  \left(\begin{array}{cc} 2 \, h_{1,2} m_{1} m_{2} + 2 \, h_{1,4} m_{1} m_{2} + 2 \, h_{1,3} m_{1} m_{3} & -2 m_{1} m_{2}( h_{1,2}  + \, h_{1,4} ) \\
-2 m_{1} m_{2}( h_{1,2}  + \, h_{1,4} ) & 2 \, h_{1,2} m_{1} m_{2} + 2 \, h_{1,4} m_{1} m_{2} + 2 \, h_{2,3} m_{2} m_{3}  \\
\end{array}\right) $$

and

$$ \mathcal{H}_2 =  \left(\begin{array}{cc}
2 m_1 (2 \, h_{1,5} m_{1} +  (h_{1,2} + h_{1,4} )m_2 + h_{1,3} m_{3})  & 2  m_{1} m_{2} ( - h_{1,2} +   h_{1,4} )\\
 -2 \, h_{1,2} m_{1} m_{2} + 2 \, h_{1,4} m_{1} m_{2} & 2 m_2 ((h_{1,2}  + h_{1,4}) m_{1} + 2 \, h_{2,4} m_{2} + h_{2,3} m_{3} )
\end{array}\right) $$

The determinant of $\mathcal{H}_1$ is 

$$ \det(\mathcal{H}_1) =4 \, {\left(( h_{1,2} + h_{1,4} )(h_{1,3} m_{1} +h_{2,3} m_{2}) + h_{1,3} h_{2,3} m_{3}\right)} m_{1} m_{2} m_{3} $$

For positive masses and $\pi/2 > \theta_1 > \theta_2 > 0$,  $\det(\mathcal{H}_1)>0$ and the trace of $\mathcal{H}_1$ is negative, so this block always has negative eigenvalues.

The submatrix $\mathcal{H}_2$ is more complicated, and it can have a positive eigenvalue for $\theta_1$ close to $\pi/2$ for some mass values.  However, using interval arithmetic we found that $\det(\mathcal{H}_2)$ is nonzero for all convex central configurations in the set $\mathcal{C}$.  This implies that all convex central configurations are minima, which in turn implies that for each positive mass vector there is at most one convex central configuration.

% need existence as well as uniqueness - ?  

\end{proof}

%%%% Yiyang 1+5 results begin

\begin{figure}[!h]
           	\centering
             \begin{tikzpicture}
             \draw (0,  0) circle (2);

32
             \draw (0, 0) -- (-35:2);
             \draw[->] (0, 0) -- (0:2);
             \draw (0, 0) -- (-80:2);
             \draw (0, 0) -- (-50:2);
             \draw (0, 0) -- (32:2);
             \draw (0, 0) -- (86:2);

             \fill (-80:2) circle (2pt) node[below] {$m_5$};
             \fill (-50:2) circle (2pt) node[right] {$m_4$};
             \fill (-35:2) circle (2pt) node[right] {$m_3$};
             \fill (32:2) circle (2pt) node[right] {$m_2$};
             \fill (86:2) circle (2pt) node[above] {$m_1$};

             \draw (.5,0) arc (0:32:.5);
             \node[] at (15:0.8) {$\theta_2$};
             \draw (1.3,0) arc (0:86:1.3);
             \node[] at (60:1.5) {$\theta_1$};
             \draw (.6,0) arc (0:-50:.6);
             
             \node[] at (-20:1.8) {$\theta_3$};
             \draw (1.6,0) arc (0:-35:1.6);
             
             \node[] at (-15:0.8) {$\theta_4$};
             \draw (1.1,0) arc (0:-80:1.1);
             \node[] at (-60:1.35) {$\theta_5$};

             \end{tikzpicture}
             \caption{Convex case for 1+5-body coorbital configuration}
             \label{fig:1}
\end{figure}
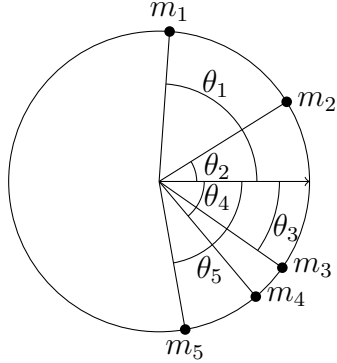

For the convex $1+5$ coorbital problem, we have the following result on a symmetry of the masses implying symmetry in the central configurations:

\begin{customthm}{6}
A convex $1+5$ coorbital central configuration ordered with $\theta_5 < \theta_4 < \theta_3 < \theta_2 < \theta_1$ and $m_1 = m_5$ and $m_2 = m_4$ must have an axis of symmetry.
\end{customthm}
\begin{proof}

Without loss of generality we can assume that  $-\frac{\pi}{2} < \theta_5 < \theta_4 < \theta_3 < \theta_2 < \theta_1 < \pi/2$.

The third row of the central configuration equations $Fm = 0$ in this case is
\begin{equation}\label{eq31}
  m_1 \sin \theta_1 (\frac{1}{r_{1,3}^3} - 1) + m_2 \sin \theta_{2,3} (\frac{1}{r_{2,3}^3} - 1) + m_2 \sin \theta_{4,3} (\frac{1}{r_{3,4}^3} - 1) + m_1\sin \theta_5 (\frac{1}{r_{3,5}^3} - 1)= 0.
\end{equation}

which is equivalent to
$$  m_1 \frac{r_{3,5}^3(1 - r_{1,3}^3) \sin \theta_{1,3} + r_{1,3}^3 (1 - r_{3,5}^3) \sin \theta_{5,3}}{r_{1,3}^3 r_{3,5}^3} +  m_2 \frac{r_{3,4}^3 (1 - r_{2,3}^3) \sin \theta_{2,3} + r_{2,3}^3 (1 - r_{3,4}^3) \sin \theta_{4,3}}{r_{2,3}^3 r_{3,4}^3} = 0.$$

The numerators of the above equation can be rewritten using Lemma \ref{c4lem}, and that lets us write the equation as

\begin{equation}\label{boxtri}
  m_2 \sin (\frac{\theta_{3,4} - \theta_{2,3}}{4}) (\triangle) + m_1 \sin (\frac{\theta_{3,5} - \theta_{1,3}}{4})(\Box) = 0,
\end{equation}
where 

\begin{align*}
\triangle & = \frac{32 \sin \frac{\theta_{2,3}}{2} \sin \frac{\theta_{3,4}}{2}}{r_{2,3}^3 r_{3,4}^3} \left[8 \sin^2 \frac{\theta_{2,3}}{2} \sin^2 \frac{\theta_{3,4}}{2} \cos(\frac{\theta_{2,3} + \theta_{3,4}}{2}) \cos ( \frac{\theta_{3,4} - \theta_{2,3}}{4} )\right . \\
& + \left . \sin (\frac{\theta_{2,3}  + \theta_{3,4}}{4})(1 + \cos \frac{\theta_{2,3}}{2} \cos \frac{\theta_{3,4}}{2}) \right]
\end{align*}

and 

\begin{align*} 
\Box & = \frac{32 \sin \frac{\theta_{1,3}}{2} \sin \frac{\theta_{5,3}}{2}}{r_{1,3}^3 r_{3,5}^3} \left[8 \sin^2 \frac{\theta_{1,3}}{2} \sin^2 \frac{\theta_{5,3}}{2} \cos(\frac{\theta_{1,3} + \theta_{5,3}}{2}) \cos (\frac{\theta_{5,3} - \theta_{1,3}}{4}) \right . \\
& \left . + \sin (\frac{\theta_{1,3} + \theta_{5,3}}{4}) (1 + \cos \frac{\theta_{1,3}}{2} \cos \frac{\theta_{5,3}}{2}) \right]
\end{align*}
both of which are positive for these convex configurations (every factor is positive for the angles under consideration). 

Next we add rows 2 and 4 from the equations $Fm = 0$, obtaining

\begin{equation}
\begin{split}\label{eq33}
   & m_1 \frac{ \sin(\theta_{1,2}) r_{4,5}^3 (1 - r_{1,2}^3) + \sin(\theta_{5,4}) r_{1,2}^3 (1 - r_{4,5}^3)} {r_{1,2}^3 r_{4,5}^3} - m_3 \frac{\sin \theta_{2,3} r_{3,4}^3 (1 - r_{2,3}^3) + \sin \theta_{4,3} r_{2,3}^3 (1 - r_{3,4}^3)}{r_{2,3}^3 r_{3,4}^3} + \\
   & m_1 \frac{\sin(\theta_{1,4}) r_{2,5}^3 (1 - r_{1,4}^3) + \sin(\theta_{5,2}) r_{1,4}^3 (1 - r_{2,5}^3)}{r_{1,4}^3 r_{2,5}^3} = 0
\end{split}
\end{equation}

We can use Lemma \ref{c4lem} to rewrite each of these three terms, and combine the first and last terms in the same way as in Theorem \ref{sym4}.  Then equation (\ref{eq33}) can be written as:

\begin{equation} \label{boxheart}
m_1 \sin (\frac{\theta_{4,5} - \theta_{1,2}}{4})(\heartsuit)  -  m_3 \sin ( \frac{\theta_{3,4}- \theta_{2,3}}{4}) (\triangle) = 0
\end{equation}
where 

\begin{align*} 
\heartsuit & = 4 \cos (\frac{\theta_{4,5} - \theta_{1,2}}{4}) \big( \cos(\frac{\theta_{1,2}+ \theta_{4,5}}{2}) + \cos(\frac{\theta_{1,4} + \theta_{2,5}}{2})\big) \\
& + \frac{1}{2 \sin^2 \frac{\theta_{1,2}}{2} \sin^2 \frac{\theta_{4,5}}{2}} \sin (\frac{\theta_{1,2} + \theta_{4,5}}{4}) (1 + \cos \frac{\theta_{1,2}}{2} \cos \frac{\theta_{4,5}}{2}) \\
& + \frac{1}{2 \sin^2 \frac{\theta_{1,4}}{2} \sin^2 \frac{\theta_{2,5}}{2}} \sin (\frac{\theta_{1,4} + \theta_{2,5}}{4}) (1 + \cos \frac{\theta_{1,4}}{2} \cos \frac{\theta_{2,5}}{2}) \\
\end{align*}
is positive (this is the same quantity shown to be positive in Theorem \ref{sym4}), and $\triangle$ is defined as above (which is also positive on these configurations).

%We also have equation \ref{boxtri}
%    
%    In the case that $\theta_{3,4}- \theta_{2,3} \geq 0$, equation \ref{boxheart} implies that $\theta_{4,5} - \theta_{1,2} \geq 0$.  Together these imply that $\theta_{3,5} - \theta_{1,3} \geq  \theta_{3,4} - \theta_{2,3} \geq 0$.   But this means each term in the left-hand side of equation \ref{boxtri} is non-negative, so $\theta_{3,5} -\theta_{1,3} = \theta_{3,4} - \theta_{2,3} = 0$. Thus, we conclude that $\theta_{1,3} = \theta_{3,5}$ and $\theta_{2,3} = \theta_{3,4}$, i.e. the 1+5-body coorbital central configuration must have an axis of symmetry.

For these convex configurations with positive masses, equation \ref{boxheart} implies that 

$$(\theta_{4,5} - \theta_{1,2}) (\theta_{3,4} - \theta_{2,3} ) \geq 0$$
while equation \ref{boxtri} implies that

$$ (\theta_{3,4} - \theta_{2,3})(\theta_{3,5} - \theta_{1,3}) \leq 0 $$

Since

$$(\theta_{3,5}-\theta_{1,3}) = (\theta_{3,4}-\theta_{2,3})+(\theta_{4,5}-\theta_{1,2})$$
this is only possible if $\theta_{3,5}-\theta_{1,3} = \theta_{3,4}-\theta_{2,3} = \theta_{4,5}-\theta_{1,2} = 0$, and the configuration is symmetric about the third point.

\end{proof}

%%%%% Yiyang 1+5 results end

\section{Conclusion and Future Work}

The $1+N$ coorbital problem seems to be rich in interesting questions.   Many of the results in this work suggest generalizations for large values of $N$; for large $N$ it seems crucial to choose coordinates which scale better than the mutual distances.   Besides their inherent interest, further results on the finiteness and enumeration of central configurations for the coorbital problem may also shed some light on the finiteness problem for the planar $N$-body problem.  

\section{Acknowledgements}
Yiyang Deng was partially supported by the Science and Technology Research Program of Chongqing Municipal Education Commission Grant No. KJQN201800815, the Research Program of CTBU Grant No. 1952040 and the Program for the Introduction of High-Level Talents of CTBU Grant No. 1856010.

\bibliography{CelMechEtc}{}
\bibliographystyle{amsplain}
\end{document}